%% file: main.tex
\documentclass[twoside]{article}

\usepackage{bbold}
\usepackage{latexsym}
\usepackage{mathtools}              % Тот же amsmath, только с некоторыми поправками
\usepackage{amssymb}                % Математические символыы
\usepackage{amstext}                % Текстовые вставки в
\usepackage{amsfonts}               % Математические шрифты
\usepackage{icomma}                 % "Умная" запятая: $0,2$ ---
\usepackage{enumitem}               % Для выравнивания itemize
\usepackage{array}                  % Таблицы и матрицы
\usepackage[matrix,arrow,curve]{xy} % Матрицы, стрелки и кривые
\usepackage[ruled,vlined,linesnumbered]{algorithm2e}            % algorithms
\usepackage{nicefrac}       % compact symbols for 1/2, etc.
\usepackage{multicol}
\usepackage{color}
\usepackage[colorlinks=true,linkcolor=red, citecolor=blue, urlcolor=blue]{hyperref}

\newtheorem{lemma}{Lemma}
\newtheorem{theorem}{Theorem}
\newtheorem{definition}{Definition}
\newtheorem{prop}{Proposition}
\newtheorem{assumption}{Assumption}
\newtheorem{corollary}{Corollary}

\SetKwInput{KwInput}{Input}                % Set the Input
\SetKwInput{KwOutput}{Output}              % set the Output

\newcommand{\E}{\mathbb{E}}
\newcommand{\argmin}{\operatornamewithlimits{argmin}}
\newcommand{\argmax}{\operatornamewithlimits{argmax}}

\newcommand{\Prob}{\operatorname{Pr}}

\newcommand{\Mirr}{\operatorname{Mirr}}

\newcommand{\one}{\mathbf{1}}
\newcommand{\bI}{\mathbf{I}}
\newcommand{\bP}{\mathbf{P}}
\newcommand{\br}{\mathbf{r}}

\newcommand{\bbR}{\mathbb{R}}

\newcommand{\Atot}{A_{\text{tot}}}
\newcommand{\tmix}{t_{\text{mix}}}
\newcommand{\St}{\mathcal{S}}
\newcommand{\Act}{\mathcal{A}}

\newcommand{\eps}{\varepsilon}
\newcommand{\poly}{\operatorname{poly}}

\usepackage[dvipsnames]{xcolor}
\newcommand{\new}[1]{{#1}}

\usepackage[accepted]{aistats2022}
\usepackage[round]{natbib}

\begin{document}

% If your paper is accepted and the title of your paper is very long,
% the style will print as headings an error message. Use the following
% command to supply a shorter title of your paper so that it can be
% used as headings.
%
%\runningtitle{I use this title instead because the last one was very long}

% If your paper is accepted and the number of authors is large, the
% style will print as headings an error message. Use the following
% command to supply a shorter version of the authors names so that
% they can be used as headings (for example, use only the surnames)
%
%\runningauthor{Surname 1, Surname 2, Surname 3, ...., Surname n}

\twocolumn[

\aistatstitle{Primal-Dual Stochastic Mirror Descent for MDPs}

\aistatsauthor{Daniil Tiapkin \And Alexander Gasnikov}

\aistatsaddress{ HSE University, Russia \And  Moscow Institute of Physics and Technology, Russia\\
HSE Univeristy, Russia} ]

\begin{abstract}
We consider the problem of learning the optimal policy for infinite-horizon Markov decision processes (MDPs). For this purpose, some variant of Stochastic Mirror Descent is proposed for convex programming problems with Lipschitz-continuous functionals. An important detail is the ability to use inexact values of functional constraints and compute the value of dual variables. We analyze this algorithm in a general case and obtain an estimate of the convergence rate that does not accumulate errors during the operation of the method. Using this algorithm, we get the first parallel algorithm for mixing average-reward MDPs with a generative model without reduction to discounted MDP. One of the main features of the presented method is low communication costs in a distributed centralized setting, even with very large networks.
\end{abstract}

\section{\uppercase{Introduction}}

We consider the following nonsmooth convex optimization problem over a simple closed convex set $Q \subseteq E$, where $E$ is a finite-dimensional normed space, with additional functional constraints
\begin{align*}%\label{eq:general_problem}
    \begin{split}
    \min_{x\in Q}& \ f(x), \\
    \text{s.t. }& g^{(l)}(x) \leq 0, \ \forall l \in [m].
    \end{split}
\end{align*}
In the context of modern large-scale optimization, the number of constraints $m$ could be huge, so we are interested in first-order algorithms in which a number of iterations does not depend on $m$. In book by \cite{nemirovski1983problem} it was noticed that the usual (sub)gradient method could handle functional constraints without additional price in terms of a total number of iterations. The proposed scheme in the simple form could be written as follows
\[
    x^{k+1} = \begin{cases}
        x^{k} - \eta \nabla f(x^k), & \max\limits_{l \in [m]} g^{(l)}(x^k) \leq \eps; \\
        x^{k} - \eta \nabla g^{(l(k))}(x^k), & \max\limits_{l \in [m]} g^{(l)}(x^k)  > \eps,
    \end{cases}
\]
where $l(k) = \argmax_{l \in [m]} g^{(l)}(x^k)$ and $\eps$ is a desired accuracy of constraint satisfaction. We emphasize that this scheme could be paralleled very efficiently: different threads or nodes of the network could handle the computation of different constraints. The ability to parallel computation is crucial for any large-scale application.

Next there were two main directions in the development of this scheme:
\begin{itemize}
    \item use of stochastic (sub)gradients;
    \item computation of dual variables for the Lagrange dual problem (primal-duality).
    %\item handle noise in constraint computations.
\end{itemize}

The first direction is essential for large-scale applications because the computation of exact gradients could be impossible or computationally heavy. The value of the second type of development highly depends on the particular application but always gives a possibility to use stopping criteria based on the duality gap.

The development of the stochastic case was initiated in the paper \citep{nemirovski2009robust} for Mirror Descent without functional constraint and developed for high-probability deviations bounds in the paper \citep{lan2012validation}. The work \citep{bayandina2018mirror} makes possible application of these results to the scheme with functional constraints.

The first step to the primal-duality of subgradient method was done in the paper \citep{nesterov2009primal-dual} but in a different direction: the author solves the dual problem using subgradient method and reconstructs the primal variables. In the work \citep{nesterov2014primal-dual} the authors propose the scheme that solves the primal conic problem using subgradient method with functional constraints and afterward computes dual variables without any computational price. This approach was generalized to arbitrary deterministic convex optimization problems in the paper \citep{bayandina2018mirror}. 

%The development of the stochastic noise in constraint computation was done by \cite{lan2020algorithms} but in the price of additional

However, a natural question appears: \textit{Is it possible to combine these two properties and propose a stochastic subgradient algorithm that computes the dual variables without additional computational price?} In this paper, we give a positive answer. Additionally, we propose a bright application that requires combining both traits with additional inexactness in constraint computation.

\paragraph{Markov Decision Process.}

We apply the proposed primal-dual stochastic Mirror Descent to the problem of \textit{mixing average-reward Markov Decision Process}. 

Markov Decision Process (MDP) is a mathematical model for the reinforcement learning (RL) problem, the rapidly developing branch of modern machine learning \citep{sutton1998rl,szepesvari2010algorl}. We consider the infinite horizon average-reward setting of this problem. For complete notations we refer to Section~\ref{sec:mixing_amdp}.

The solution to the average-reward MDP (AMDP) with $S$ states and $\Atot$ state-action pairs could be described through the linear program that obtained from Bellman equations \citep{bertsekas2005dynamic}
\begin{align*}%\label{eq:amdp_lp_vanila}
    \begin{split}
        \min_{\bar v, h} &\  \bar v \\
        \text{s.t. } &\bar v \one + (\hat\bI- \bP) h - \br \geq 0,
    \end{split}
\end{align*}
where \(\hat \bI_{(i,a_i),j} = \bI_{i,j}, \bP \in \bbR^{\Atot \times S}\) is a transition probability matrix,  $\br \in \bbR^{\Atot}$ is a vector of rewards for state-action pairs, $\bar v$ is an average reward value, and $h$ is a bias vector. For another setting of discounted MDP we refer to the line of the previous work \citep{azar2012sample,sidford2018near,sidford2018variance,agarwal2020model,li2020breaking} and references within.

Let us list important properties of the problem \eqref{eq:amdp_lp_vanila}: 1) it has a huge number of constraints $\Atot$; 2) it is impossible to compute neither constraints nor its gradient since the transition probability matrix $\bP$ is usually unknown; 3) a policy that corresponds to the optimal average reward $\bar v$ can be computed from the solution to the dual problem of \eqref{eq:amdp_lp_vanila}. 

To handle the second problem we consider solving AMDP with \textit{generative model} or \textit{sampler} \citep{azar2012sample,jin2020efficiently,agarwal2020model}: a stochastic oracle generates a transition state from a given state-action pair according to probabilities \(\mathcal{P}\). This assumption makes applying of the proposed stochastic primal-dual Mirror Descent possible to the problem \eqref{eq:amdp_lp_vanila} after a suitable approximation of constraint functions. We underline that it is required to use a combination of \textit{all} properties of our algorithm to solve this problem. Additionally, we notice that the work \citep{lan2020algorithms} can handle stochastic constraints without additional approximation of constraint functions but at the price of much worse complexity $\sim \eps^{-4}$.

However, it is impossible to obtain rates for AMDP solving without additional assumptions. We consider the \textit{mixing} assumption that the Markov chain that corresponds to the choice of any policy converges to the stationary distribution sufficiently fast. In papers \citep{wang2017primaldual,jin2020efficiently} authors showed that under mixing assumption, the search space for the bias vector $h$ could be bounded using mixing time $\tmix$. A similar assumption was studied in the paper \citep{kearns2002near} and an alternative one of communicating MDP with a finite diameter in works \citep{bartlett2009regal,jaksch2010near,agrawal2017optimistic}.

To the authors' best knowledge, there are only three works that consider the infinite-horizon mixing AMDP with a generative model -- \citep{wang2017primaldual,jin2020efficiently,jin2021tight}. In the first two papers the same general convex optimization algorithm was used -- Stochastic Mirror Descent for saddle-point problems, and both of the presented algorithms are not designed for parallel computations. In the last paper authors perform reduction to the discounted problem, and the price of this reduction is sample complexity of order \(O(\eps^{-3})\). In this paper, we present the first \textit{parallel} algorithm for this problem without reductions to discounted MDP with very low communication costs. The parallelism gives a possibility to handle very large setups of MDPs that cannot be stored in the memory of one machine. In the case of simultaneous working of \(\Atot\) workers, our algorithm works in \(\tilde O(t^2_{\text{mix}} \vert \St\vert \eps^{-2})\) real time and outperforms approach of \cite{jin2020efficiently} which works in \(\tilde O(t^2_{\text{mix}} \Atot \eps^{-2})\).

\paragraph{Our contribution.} 
\begin{itemize}
    \item The first (sub)gradient-based algorithm for optimization with functional constraints that 1) allows the use of stochastic gradients, 2) computes the value of Lagrange dual variables, 3) allows inexact computation of constraints at the same time.
    \item The first parallel algorithm for solving mixing average-reward MDP without reduction to the discounted problem. Additionally, this algorithm has very low communication costs and thus could work on very large centralized networks effectively.
\end{itemize}

\paragraph{Paper organization.} Section~\ref{sec:smd} is devoted to the proposed primal-dual stochastic variant of the Mirror Descent algorithm. All proofs are presented in the supplementary material. Section~\ref{sec:mixing_amdp} describes how to apply the results of Section~\ref{sec:smd} to the mixing AMDP problem. Finally, Section~\ref{sec:exp} contains a numerical comparison to the approach described in the paper of \cite{jin2020efficiently}.

\paragraph{Notation.} For a matrix \(A \in \mathbb{R}^{n \times m}\) we define its \(i\)-th row as  \(A_{(i)}\). We denote by \(\one = (1,\ldots,1)^\top\) the vector filled with ones. By \(e_{i}\) we define a standard basis vector. Also we define \(\Delta^n = \{ x \in \mathbb{R}^n \mid \forall i: x_i \geq 0, \sum_{i=1}^n x_i = 1\}\) and \(\mathbb{B}^n_c = [-c,c]^n\). \(\bI\) is an identity matrix of size deducible from the context. Inner product \(\langle \cdot, \cdot \rangle \colon E^* \times E \to \mathbb{R}\) is defined on pairs of vectors from dual and primal spaces. In the case of Euclidean spaces, it coincides with the standard inner product. For a normed space $(E, \Vert \cdot \Vert)$ we define a dual norm on a space $E^*$ as follows: $\Vert v \Vert_* = \sup_{x: \Vert x \Vert = 1} \langle v, x \rangle$. Also we define \([m] = \{1,\ldots,m\}\). By \(\mathbb{I}\{A\}\) we define an indicator of a set \(A\). For $i \in [m]$ define $e_i \in \mathbb{R}^m: e_i[j] = \mathbb{I}\{i=j\}$. 
%%%
\section{\uppercase{Primal-Dual Stochastic Mirror Descent}}\label{sec:smd}

In this section we develop techniques of \cite{bayandina2018mirror}. Firstly, we introduce a basic notation that will be used further. Then we propose a new algorithm for the constrained convex stochastic optimization problem in the model of inexact computation of constraint functions. Finally, we prove convergence of the algorithm in terms of the duality gap between primal and (Lagrange) dual problems. The last part is crucial for an application on average-reward MDPs.

\subsection{Notation}

We consider the constrained convex optimization problem over a convex compact set \(Q \subseteq E\), where \(E\) is a finite-dimensional normed space
\begin{align}\label{eq:general_problem}
    \begin{split}
    \min_{x\in Q}& \ f(x), \\
    \text{s.t. }& g^{(l)}(x) \leq 0, \ \forall l \in [m],
    \end{split}
\end{align}
and \(f \colon Q \to \mathbb{R}, g^{(l)} \colon Q \to \mathbb{R}\) are convex functions. We assume that subgradients of these functions exist for each \(x\in Q\) for simplicity. We call \(\nabla f(x), \nabla g^{(l)}(x)\) any subgradients of corresponding functions. However, in our algorithm we have an access only to stochastic subgradient oracles \(\nabla f(x, \xi), \nabla g^{(l)}(x, \xi_{(l)})\). 

Now we are going to introduce definitions that will be useful in the algorithm description.
\begin{definition}\label{df:stochastic_delta_approx}
    Function \(h_\delta \colon Q \to \mathbb{R}\) is called a \(\delta\)-approximation of \(h \colon Q \to \mathbb{R}\) if \(\vert h_\delta(x) - h(x) \vert \leq \delta\) for all \(x \in Q\). 
\end{definition}
For our problem, we suggest that we have an oracle not for computation of constraints \(g^{(l)}\) but their \(\delta\)-approximations \(g^{(l)}_\delta\). It is the important difference with the setup of \cite{bayandina2018mirror}: in our assumption we cannot consider only one constraint of form \(g(x) = \max_{l\in[m]} g^{(l)}(x)\) because it does not seem possible to compute a subgradient of \(g\) given subgradients of \(g^{(l)}\) when index on which the maximal value attains is unknown.

Now we are ready to list all our assumptions on problem setup \eqref{eq:general_problem} for our algorithm:
\begin{itemize}
    \item[\textbf{(A1)}] \(f\) and all \(g^{(l)}\) are Lipschitz continuous with constant \(M\) for objective function and for all constraints;
    \item[\textbf{(A2)}] Stochastic subgradients are unbiased: \(\E \nabla f(x,\xi) = \nabla f(x), \ \E \nabla g^{(l)}(x,\xi_{(l)}) = \nabla g^{(l)}(x) \);
    \item[\textbf{(A3)}] Stochastic subgradiens are bounded: \(\Vert \nabla f(x,\xi) \Vert_* \leq M, \Vert \nabla g^{(l)}(x,\xi_{(l)}) \Vert_* \leq M\) a.s.;
\end{itemize}

Since our algorithm is Mirror-Descent based, the next step is to define the proximal step and basic properties of Mirror Descent. 

Firstly, we define a \textit{prox-function} \(d \colon Q \to \mathbb{R}\) as a continuous 1-strongly convex function \(d\) with respect to the norm \(\Vert \cdot \Vert\) on \(E\) that admits continous selection of subgradients \(\nabla d(x)\) where they exist. \textit{Bregman divergence} that corresponds to a prox-function \(d\) is a function \(V(x,y) = d(y) - d(x) - \langle \nabla d(x), y-x \rangle\). 

Given vectors \(x \in E\) and \(v \in E^*\), the mirror step is defined as
\[
    x^+ = \Mirr(x, v) = \argmin_{y \in X} \left\{\langle v, y \rangle + V(x,y) \right\}.
\]
We assume that the mirror step can be easily computed.

\subsection{Primal Problem}

In this subsection, we consider problem \eqref{eq:general_problem} in terms of convergence of the objective function in confidence region. Formally speaking, a vector \(\hat x\) is called an \((\eps_f, \eps_g, \sigma)\)-solution to the primal problem \eqref{eq:general_problem}, if
\begin{align}\label{eq:primal_criteria}
    \begin{split}
        f(\hat x) - f(x^*) &\leq \eps_f, \\
         g^{(l)}(\hat x) &\leq \eps_g \ \forall l \in \{1,\ldots,m\} \ \text{w.p.}\geq 1 - \sigma,
    \end{split}
\end{align}
where \(x^*\) is a true minimizer of the problem \eqref{eq:general_problem}. We assume that an algorithm have access only to stochastic subgradient oracles of functions \(f, g^{(l)}\) and to \(\delta\)-approximations \(g^{(l)}_\delta\) of constraint functions.

\begin{algorithm} \label{alg:md}
    \caption{Stochastic Mirror Descent with noisy constraints}
    \KwInput{accuracy \(\eps > 0\), number of steps \(N\), stepsize \(\eta=\varepsilon/M^2\)}
    %\(\eta = \eps/M^2\)\;
    \(x^0 = \argmin_{x \in Q} d(x)\)\;
    \(I = \emptyset, J = \emptyset\)\;
    \For{ \(k=0,1,2,\ldots,N-1\) }{
        \If{\(g_\delta^{(l)}(x^k) \leq \eps + \delta \ \forall l \in [m]\)}{
            \(x^{k+1} = \Mirr(x^k, \eta \nabla f(x^k, \xi^k))\) \tcp*{"productive" steps}  
            Add \(k\) to \(I\)
        }
        \Else {
            \(l(k) = \argmax\limits_{l \in [m]} g_\delta^{(l)}(x^k)\)\;
            \(x^{k+1} = \Mirr(x^k, \eta \nabla g^{(l(k))}(x^k, \xi^k_{(l(k))}))\) \tcp*{"non-productive" steps}
            Add \(k\) to \(J\);
        }
    }
    \KwRet \(\hat x = \frac{1}{\vert I \vert } \sum_{k \in I} x^k\);
\end{algorithm}

    Denote \(\hat \nabla_k f = \nabla f(x^k, \xi^k), \nabla_k f = \nabla f(x^k)\) and \(\hat \nabla_k g^{(l)} = \nabla g^{(l)}(x^k, \xi^k_{(l)}), \nabla_k g^{(l)} = \nabla g^{(l)}(x^k)\). Additionally, we define  \(\Theta_0^2 = d(x^*) - d(x^0)\). In these terms we could provide the main theorem.
    \begin{theorem}\label{th:primal_convergence}
        Algorithm \ref{alg:md} \new{with a constant stepsize \(\eta = \varepsilon/M^2\)} outputs \((\eps, \eps + 2\delta, \sigma)\)-solution for any \(\eps > 0, \sigma \in (0,1), \delta \geq 0\) in sense of \eqref{eq:primal_criteria} after     \[
            N \geq N_0 =  \frac{280 \cdot \Theta_0^2 M^2 \log(1/\sigma)}{\eps^2} .
        \]
\end{theorem}
The proof of this theorem is given in supplementary material. 

% We can give simple description of \(N_0\) in terms of O-notation:
% \[
%     N_0 = O\left(\frac{\Theta_0^2 M^2 \log(1/\sigma)}{\eps^2}\right).
% \]

\textit{Remark 1.} Notice that from theoretical point of view selection of the maximum in line 9 of Algorithm \ref{alg:md} could not be avoided. However, in case of non-stochastic constraint computation one of used heuristic is to set \(l(k)\) to \textit{any} index of violated constraint and we suggest that such heuristic could work in \new{the} noisy setting too.

\textit{Remark 2.} Notice that a quantity \(\Theta_0^2\) is not used in the pseudocode of Algorithm \ref{alg:md}. Thus, it is possible to use another initial \(x^0\) and have a warm start such that our algorithm converges faster. This warm start could be obtained by running an algorithm several times with a decreasing value of \(\eps\). 

\new{\textit{Remark 3.} We used a constant stepsize for the sake of simplicity. It is possible to adapt techniques of \cite{bayandina2018mirror,stonyakin2019adaptive} to use adaptive stepsizes that does not rely on knowledge of Lipschitz constant $M$.}

\subsection{Primal-Dual Convergence}

In this subsection, we extend properties of the previous algorithm and prove its primal-duality. First of all, let us define the (Lagrange) dual optimization problem associated with the problem \eqref{eq:general_problem}
\begin{align}\label{eq:dual_problem}
    \max_{\lambda \in \mathbb{R}^m_+} \left\{ \phi(\lambda) := \min_{x \in Q}\{ f(x) + \sum_{i=1}^n \lambda_l g^{(l)}(x) \} \right\}.
\end{align}
Call \(\lambda^*\) a solution to this dual problem (if it exists). We refer to \citep{boyd2004convex} for an additional background and examples.

It is well-known that for any \(x \in Q: g^{(l)}(x) \leq 0\ \forall l \in \{1,\ldots,m\}\) and \(\lambda \in \mathbb{R}^m_+\) \textit{the weak duality} holds:
\(
    \Delta(x,\lambda) = f(x) - \phi(\lambda) \geq 0,
\)
where \(\Delta\) is so-called \textit{the duality gap}. We assume that for our primal problem \eqref{eq:general_problem} the Slater's condition holds, i.e. \(\exists x \in Q: \forall l \in \{1,\ldots, m\}: g^{(l)}(x) < 0\). It implies that the dual problem has a solution and there is \textit{the strong duality}: \(\Delta(x^*, \lambda^*) = 0\) for any \(x^*\) and \(\lambda^*\) are solutions to the primal and the dual problems respectively.

It gives us a natural way to measure a quality of the pair \((\hat x, \hat \lambda)\) by the value of the duality gap \(\Delta(\hat x, \hat \lambda)\). Let us call the pair \((\hat x, \hat \lambda)\) a primal-dual \((\eps_\Delta, \eps_g, \sigma)\)-solution to \eqref{eq:general_problem} if the following holds with probability at least \(1-\sigma\)
\begin{align}\label{eq:primal-dual_criteria}
    \begin{split}
        \Delta(\hat x, \hat \lambda) &\leq \eps_\Delta,\\
        g^{(l)}(\hat x) &\leq \eps_g \ \ \ \forall l \in [m].
    \end{split}
\end{align}

Notice that since \(\hat x\) is not a feasible solution to the primal problem \eqref{eq:general_problem}, we do not have the weak duality inequality \(\Delta(\hat x, \hat \lambda) \geq 0\). However, the value of duality gap could be controlled from below because of controlled unfeasibility.

The most powerful property of Algorithm \ref{alg:md} is a possibility to generate a pair of primal-dual solutions in sense of \eqref{eq:primal-dual_criteria}: we could control the value of the duality gap without the explicit access to the constraint functions. 

\new{Following \citep{nesterov2014primal-dual,bayandina2018mirror}, we choose the following \(\hat \lambda \in \mathbb{R}^m_{+}\) as an estimate of dual variables}
\begin{equation}\label{eq:dual_variables}
    \hat \lambda_l = \frac{1}{\vert I \vert }\sum_{k \in J} \mathbb{I}\{l = l(k)\}
\end{equation} 
in terms of Algorithm \ref{alg:md}. Additionally, we define useful constant \(\overline{\Theta}_0^2 = \sup_{y \in Q} (d(y) - d(x^0))\).

Using \(\hat \lambda\), we could provide primal-dual properties of Algorithm \ref{alg:md}.
\begin{theorem}\label{th:primal-dual_convergence}
    Let us choose \(\hat \lambda \in \mathbb{R}^m_{+}\) as defined in \eqref{eq:dual_variables} and \(\hat x\) is an output of Algorithm \ref{alg:md} \new{with a constant stepsize $\eta = \varepsilon/M^2$}. Then the pair \((\hat x, \hat \lambda)\) is an \((\eps, \eps + 2\delta, \sigma)\)-solution in sense of \eqref{eq:primal-dual_criteria} for any \(\eps > 0, \delta \geq 0, \sigma \in (0, 1/2)\) after \[
        N \geq N_0' = \frac{128 \overline{\Theta}_0^2 M^2 (17 \log(2/\sigma) + 2 \kappa(E^*))}{\eps^2},
    \]
    where \(\kappa(E^*)\) is a constant of Nemirovski's inequality \citep{boucheron2013concentration} for the dual space.
\end{theorem}
\textit{Remark 1.} If \(E\) has a finite dimension \(d\), then we always have \(\kappa(E^*) \leq d\). Additionally, if \(E\) is endowed with \(\ell_p\) norm, then \(E^*\) is endowed with \(\ell_q\) norm, where \(1/p + 1/q = 1\), and there is a more precise bound, according to \citep{lutz2010nemirovski}
\[
    \kappa(E^*) \leq K\left(\frac{p}{p-1}, d\right) = \begin{cases}
        d^{\frac{2}{p}-1}, & p \in [1,2]\\
        d^{1-\frac{2}{p}}, & p \in (2,+\infty] \\
    \end{cases}
\]
In particular, if \(E\) has \(\ell_2\) norm, \(\kappa(E^*) = 1\). For \(p\in[2,+\infty]\) this bound is tight, however, in the case \(p\in[1,2]\) and \(d\geq 3\) it could be improved  \citep{boucheron2013concentration} to, for instance, a logarithmic bound \(\kappa(E^*) \leq 2e \log(d)- e\), that could be useful in the case of \(\ell_1\)-norm and an entropy prox-function.

We can write bound on \(N_0'\) using O-notation as follows
\[
    N = O\left(\frac{\overline{\Theta}_0^2 M^2 (\log(1/\sigma) + \kappa(E^*))}{\eps^2}\right).
\]
The only difference between primal and dual case is connected to the constant in Nemirovski's inequality.

\textit{Remark 2.} As in the primal case, in the complexity bounds we have a constant \(\overline{\Theta}_0^2\) that does not appear in the algorithm description. This fact gives us a chance to work much better in practice than using worst-case constant. The same situation with constant \(\kappa(E^*)\).

\section{\uppercase{Mixing AMDP}}\label{sec:mixing_amdp}

In this section, we discuss the application of the developed algorithm to the problem of approximate solving mixing average-reward Markov Decision Processes (MDP). Firstly, we propose basic definitions connected to MDPs. Next, we discuss some technical nuances that will appear in the algorithm and, finally, we outline the complete algorithm and its parallel implementation.

\subsection{Markov Decision Process}\label{subsec:mdp}

An instance of \textit{MDP} is a tuple \(\mathcal{M} = (\St,\Act,\mathcal{P}, \mathbf{r})\), where \(\St\) is a \textit{finite} set of states; \(\Act = \bigsqcup_{i \in \St} \Act_i\) is a finite state of actions, each set \(\Act_i\) contains actions from the state \(i\). \(\mathcal{P}\) is the collection of state-to-state transition probabilities given actions: \(\mathcal{P} = \{p_{ij}(a_i) \mid i,j \in \St, a_i \in \Act\}\) where \(p_{ij}(a_i)\) is a probability of transition from a state \(i\) to a state \(j\) given an action \(a_i\). Also, we define \(\mathbf{r} \in [0,1]^{\vert \Act \vert}\) as the state-action reward vector, \(r_{i,a_i}\) is the instant reward received when taking the action \(a_i\) at the state \(i \in \St\). For consistency of notation with work of \cite{jin2020efficiently}, let \((i,a_i) \in \mathcal{A}\) denote an action \(a_i\) at a state \(i\). \(\Atot = \vert \Act \vert = \sum_{i \in \St} \vert \Act_i\vert \) denotes the total number of state-action pairs. Also we denote by \(\bP \) the action-state transition probability matrix of size \(\Atot \times \vert \St \vert\), where \(\bP_{(i,a_i),j} = p_{ij}(a_i)\) in terms of \(\mathcal{P}\).

The goal is to find a stationary (randomized) policy that specifies actions to choose in the fixed state. Formally, a policy \(\pi\) is a block vector such that \(i\)-th block corresponds to a probability distribution over \(\mathcal{A}_i\). Define as \(\bP^\pi, \br^\pi\) the transition matrix and the cost vector under the fixed policy \(\pi\).

Now we are going to define optimality of the policy. In this paper we consider the infinite-horizon \textit{average-reward} MDP with the following objective to maximize
\[
    \bar v^\pi = \lim_{T \to \infty } \E_\pi\left[ \frac{1}{T} \sum_{t=1}^T \br_{i_t, a_{t}} \mid i_1 \sim \mathbf{q} \right].
\]
Here \(\{i_1,a_1,\ldots,i_t,a_t\}\) are state-actions transitions generated by MDP under a policy \(\pi\), \(\mathbf{q}\) is an initial distribution, and an expectation \(\E_\pi[\cdot]\) is taken over trajectories. In our case we interested in the case then the Markov chain generated by an AMDP under a fixed policy \(\pi\) has a unique stationary distribution \(\nu^\pi: \nu^\pi \cdot \bP^\pi = \nu^\pi\). Notice that in this case the value of \(\bar v^\pi\) does not depend on an initial distribution \(\mathbf{q}\). Then the objective simplifies a lot:
\[
    \bar v^\pi = \langle \nu^\pi, \br^\pi \rangle.
\]

Next, we define Bellman equations for an AMDP \citep{bertsekas2005dynamic}: \(\bar v^*\) is the optimal average reward if and only if there exists a vector \(h^* \in \mathbb{R}^{\vert \St \vert}\) satisfying the following
\[
    \bar v^* + h_i^* = \max_{a_i \in \Act_i}\left\{ \sum_{j \in \St} p_{ij}(a_i) h_j^* + r_{i, a_i}  \right\}, \forall i \in \St.
\]

We focus on study of the primal LP which solution is equivalent to the solution to Bellman equation
\begin{align}\label{eq:amdp_lp_vanila}
    \begin{split}
        \min_{\bar v, h} &\  \bar v \\
        \text{s.t. } &\bar v \one + (\hat\bI- \bP) h - \br \geq 0,
    \end{split}
\end{align}
where \(\hat \bI_{(i,a_i),j} = \bI_{i,j}\).

However, without additional assumptions it is hard to analyze problem. Following \citep{jin2020efficiently,wang2017primaldual}, we introduce one important assumption on an AMDP instance.
\begin{assumption}[Mixing AMDP]
    We call an AMDP instance mixing if its so-called mixing time (defined below) is bounded
    \[
        t_{\text{mix}}:= \max_{\pi}\left[ \argmin_{t \geq 1} \left\{ \max_{\mathbf{q}} \Vert ({\bP^\pi}^\top)^t \mathbf{q} - \nu^\pi \Vert_1 \leq \frac{1}{2} \right\}\right].
    \]
\end{assumption}

The most powerful corollary of this result is a possibility to make the search space a compact convex set. Formally speaking, in \citep{jin2020efficiently} it was proven that the search space for primal variables could be reduced to \(\mathcal{X} = [0,1] \times \mathbb{B}^{\vert \St \vert}_{2R} = [0,1] \times [-2R,2R]^{\vert \St \vert}\), where \(R = 2 t_{\text{mix}}\). We choose bounds of size \(2R\) instead of \(R\) because of the same reasons as \citep{jin2020efficiently} that will be described in Section \ref{sec:rounding}.

Overall, we have a (linear) optimization problem on a compact with a large number of constraints
\begin{align}\label{eq:amdp_lp}
    \begin{split}
        \min_{\bar v, h \in \mathcal{X}} &\  \bar v \\
        \text{s.t. } &\bar v \one + (\hat\bI- \bP) h - \br \geq 0.    
    \end{split}
\end{align}

If we knew the matrix \(\bP\), we could apply Stochastic Mirror Descent with constraints \citep{bayandina2018mirror} and get an approximate solution to this LP problem. However, it is not the case: we have only sampling access to the transition probability matrix. Another problem we face is computing an optimal policy by an approximate solution to this linear program. It is known that there is a strong connection between the optimal policy and the dual LP  but not the primal one. Therefore, we use primal-duality of Algorithm \ref{alg:md} and construct a policy using dual variables.

\subsection{Preprocessing}\label{sec:preprocessing}

In this subsection we aim to describe complexity of the preprocessing connected to the estimate of the transition probability matrix. We take the estimate of the form
\[
    \widetilde \bP_{(i, a_i)} = \frac{1}{n} \sum_{j=1}^n e_{X_j},
\]
where \(X_j\) are sampled from categorical distribution \(\bP_{(i, a_i)}\). Choosing appropriate \(N\) and compute these quantities for each state-action pair in parallel, we obtain the following proposition.

\begin{prop}\label{lm:preprocessing}
    For each \(\delta',\sigma' > 0\), the estimate \(\widetilde \bP\) of \(\bP\), such that for each \(a \in \Act, h \in \mathbb{B}^{\vert \St \vert}_{2R}: \vert \langle \bP_{(a)} - \widetilde \bP_{(a)}, h \rangle \vert \leq \delta'\) with probability at least \(1- \sigma'\), could be computed in \(O\left(t_{\text{mix}}^2 \Atot \cdot \frac{\vert \St \vert + \log(\Atot/\sigma')}{\delta'^2} \right)\) total samples, \(O(1)\) parallel depth and \(O\left( t_{\text{mix}}^2\frac{\vert \St \vert  + \log(\Atot/\sigma')}{\delta'^2} \right)\) samples proceed by each single node. { In the case of \(m \leq \Atot\) available workers, it works in \( O\left( \frac{\Atot}{m} \cdot t_{\text{mix}}^2\frac{\vert \St \vert  + \log(\Atot/\sigma')}{\delta'^2} \right)\) real time.}
\end{prop}
The proof of this proposition is given in supplementary.

\subsection{Rounding to Optimal Policy}\label{sec:rounding}

In this subsection, we prove the result that give us a possibility to obtain an approximate optimal policy from the dual variables produced by \eqref{alg:md}.

\begin{prop}\label{lm:rounding}
    Suppose that primal \((\bar v^\eps, h^\eps)\) and dual \(\mu^\eps\) variables are \((\eps_f, \eps_g, \sigma)\)~-approximate solution to \eqref{eq:amdp_lp} in terms of (expected) primal-dual convergence \eqref{eq:primal-dual_criteria}. Define the policy \(\pi\): \(\pi_{i,a_i} = \frac{\mu_{i,a_i} }{\lambda_i}\), where \(\lambda \in \mathbb{R}_+^{\vert \St \vert}\) is defined as \(\lambda_i = \sum_{a_i \in \mathcal{A}_i} \mu_{i,a_i}\). Then \(\pi\) is an \(4(\eps_f + \eps_g)\)-optimal policy with probability at least \(1-\sigma\).
\end{prop}
The proof is given in supplementary material, and it is very similar to the proof of \citep{jin2020efficiently}. There are two differences. Firstly, we have guarantees on our primal-dual solution, not the solution to a saddle-point problem, and this slightly changes the structure of the proof. Secondly, we have high-probability bounds instead of bounds in expectation. 

\subsection{Parallel Algorithm}

In this subsection, we describe a final algorithm to approximate solving AMDP in parallel. In our setup, each single node of a centralized network corresponds to a single state-action pair.

First of all, we describe linear program corresponding to AMDP \eqref{eq:amdp_lp} in terms of \eqref{eq:general_problem},
\begin{align}
    \begin{split}\label{eq:amdp_reform}
        \min_{\bar v, h \in \mathcal{X}} &\  f(\bar v, h) = \bar v, \\
        \text{s.t. } &g^{(i, a_i)}(\bar v, h) = r_{(i, a_i)} - \bar v\\
        &\qquad + (\langle \mathbf{P}_{(i, a_i)}, h \rangle - h_i) \leq 0  \quad \forall (i,a_i) \in \mathcal{A}.
    \end{split}
\end{align}
where \(\mathcal{X} = [0,1] \times \mathbb{B}_{4\tmix}^{\vert \St \vert}\). We set standard Euclidean prox-structure on \(\mathcal{X}\): \(\ell_2\)-norm \(\Vert \cdot \Vert_2\) and a prox-function \(d(x) = \frac{1}{2}\Vert \cdot \Vert_2^2\). In these terms, we have the following constants: \(M = 2\) and \(\overline{\Theta}_0^2 = (4R)^2 \vert S \vert + 1 = O(\tmix^2 \vert S \vert)\). 

\new{However, in our case we cannot compute constraints since there is no access to the true model. To overcome this, we run preprocessing described in Section \ref{sec:preprocessing} to derive approximate model $\widetilde \bP$ and obtain \(\delta\)-approximation of constraints \(g^{(i,a_i)}_{\delta}\)
\begin{equation}\label{eq:g_approx}
    g^{(i,a_i)}_{\delta}(\bar v,h) = r_{(i,a_i)} - \bar v + (\langle \widetilde \bP_{(i,a_i)}, h \rangle - h_i).
\end{equation}

Additionally, we are going to use stochastic subgradients for \(g\) by using samples of next state \(s \sim \mathbf{P}_{(i,a_i)}\):
\begin{small}
\begin{align}\label{eq:stochastic_gradients_amdp}
     &\nabla_{\bar v} g^{(i,a_i)}(\bar v, h) = -1, 
     &\hat \nabla_{h} g^{(i,a_i)}(\bar v, h, s) = e_{s} - e_i,
 \end{align}
 \end{small}
 with constant \(M = 2\). In this case, we can use Algorithm \ref{alg:md} with the following update rules for primal variables $\bar v^k$ and $h^k$. For productive steps $(k \in I)$
 \begin{equation}\label{eq:primal_update_prod}
     \bar v^{k+1} = \bar v^{k} - \eta, \quad h^{k+1} = h^{k},
\end{equation}
and for non-productive steps $(k \in J)$
\begin{equation}\label{eq:primal_update_nonprod}
     \bar v^{k+1} = \bar v^{k} + \eta, \quad h^{k+1} = h^{k} - \eta(e_s - e_i),
\end{equation}
 where $(i, a_i)$ is a state-action pair where $\max_{a \in \Act} g_\delta^{a}(\bar v^k, h^k)$ attains, and $s \sim \bP_{(i,a_i)}$ is a transition sample. We note that in this form this algorithm is sequential.
 
 To design a parallel version of Mirror Descent that presented in Algorithm~\ref{alg:amdp_solving}, we are going to use a separate node for each state-action pair $(j,a_j)$. The aim of node  corresponds to state-action pair $(j,a_j)$ is  1) compute and update value of $c^k_{(j,a_j)} := g_\delta^{(j,a_j)}$ and 2) sample transitions of MDP. 

To compute $c^k_{(j,a_j)}$ faster than in $O(\vert \St \vert)$ operations, we note that updates of $\bar v^k$ and $h^k$ are sparse and we can store the previous value $c^{k-1}_{(j, a_j)}$.

In the case of productive steps $(k \in I)$ we update as follows:
\begin{equation}\label{eq:c_k_update_prod}
    c^k_{(j, a_j)} = c^{k-1}_{(j,a_j)} + \eta.
\end{equation}
This update rule is correct by the update rule \eqref{eq:primal_update_prod} and equation \eqref{eq:g_approx}.

In the case of non-productive steps $(k \in J)$ we have a more complicated update
\begin{equation}\label{eq:c_k_update_nonprod}
    \begin{split}
        c^k_{(j,a_j)} &= c^{k-1}_{(j,a_j)} - \eta(1 + \widetilde{\bP}_{(j,a_j),s } - \widetilde{\bP}_{(j,a_j), i} \\
        &+ \mathbb{I}\{j = i\} - \mathbb{I}\{j = s\}).
    \end{split}
\end{equation} 
The correctness of the procedure is guaranteed by the update rule \eqref{eq:primal_update_nonprod} and definition \eqref{eq:g_approx}. Update rules \eqref{eq:c_k_update_prod} and \eqref{eq:c_k_update_nonprod} gives us an opportunity to update constraints in $O(1)$ time per each node.

\begin{theorem}
    Let \(\eps > 0\) and \(\sigma \in (0, 1/2)\). The policy \(\hat \pi\) generated by Algorithm \ref{alg:amdp_solving} \new{with a constant stepsize $\eta=\eps/64$} and preprocessing described in Section \ref{sec:preprocessing} performed with parameters \(\delta' = \eps/16, \sigma' = \sigma/2\) is an \(\eps\)-approximate optimal policy with probability at least \(1 - \sigma\) if
    \[
        N =  O\left( \frac{\tmix^2  \vert \St \vert \log(1/\sigma)}{\eps^2} \right).
    \]
    The described algorithm has \(O(1)\) parallel depth, \( O(\Atot \cdot N)\) sample and running time complexity.  In the case of \(m \leq \Atot\) available workers, algorithm works in \(O\left(\frac{\Atot}{m} \cdot t_{\text{mix}}^2 \vert \St \vert \log(1/\sigma) \cdot \eps^{-2}\right)\) real time. The full description is presented in Algorithm~\ref{alg:amdp_solving}.
\end{theorem}
}

\begin{algorithm}
    \caption{Parallel average-reward MDP}
    \label{alg:amdp_solving}
    \SetAlgoLined\DontPrintSemicolon
    \KwInput{accuracy \(\eps > 0\), number of steps \(N\), 
    approximation accuracy $\delta = \eps/16$,
    Mirror Descent accuracy $\tilde\eps = \eps/16$, stepsize $\eta=\tilde\eps/4$, confidence level $\sigma$.}
    % Set Function Names
    \SetKwFunction{FHead}{MirrorDescent}
    \SetKwFunction{FWorker}{WorkerNode}
     
      \SetKwProg{Fn}{Procedure}{:}{}
      \Fn{\FHead{}}{
            \(\bar v^0 = 0; h^0 = 0\);\;
            \For{\(k = 0,1,2\ldots, N-1\)} {
                Send "Check constraints" to all nodes;\;
                Receive \(g^{(i, a_i)}_{\delta}(\bar v^k, h^k) \ \forall (i, a_i) \in \mathcal{A}\);\;
                \If{\(\max\limits_{a \in \Act} g^{(i, a_i)}_{\delta}(\bar v^k, h^k) \leq \tilde\eps + \delta\)} {
                    \(\bar v^{k+1} = \bar v^{k} - \eta\),\ \ \(h^{k+1} = h^k\);\;
                    Send "Productive step" to all nodes;\;
                    Add \(k\) to \(I\);\;
                } \Else {
                    \((i,a_i)_k = \argmax_{a \in \Act} g_{\delta}^{(i,a_i)}(\bar v^k, h^k)\);\;
                    Send "Sample" to node \((i,a_i)_k\);\;
                    Receive state \(s\);\;
                    \(\bar v^{k+1} = \bar v^{k} + \eta\), \ \ \(h^{k+1} = h^k - \eta \cdot (e_{s} - e_i)\);\;
                    Send "Non-productive step", \(i\) and \(s\) to all nodes;\;
                    Add \(k\) to \(J\);\;
                }
            }
            \(\hat \mu_{(i,a_i)} = \frac{1}{\vert I \vert} \sum_{k \in J} \mathbb{I}\{(i,a_i)_k = (i,a_i)\}\);\;
            \KwRet \(\pi_{i,a_i} = \mu_{i,a_i} / \left(\sum_{a_i \in \Act_i} \mu_{i, a_i}\right)\);
     }
    \SetKwProg{Fn}{Procedure}{:}{}
    \Fn{\FWorker{\(j, a_j\)}}{
            Compute \(\widetilde \bP_{(j, a_j)}\) with precision $\delta$ and confidence level $\sigma/2$;\;
            \(c^0_{(j, a_j)} = r_{j,a_j}\),\(k = 0\);\;
            \While{is not finished} {
                Wait message;\;
                \If{"Compute constraints"}{
                    Send \(c^k_{(j,a_j)}\) as \(g^{(j,a_j)}_{\delta}(\bar v^k, h^k)\);\;
                }  \ElseIf{"Sample"} {
                    Sample \(s \sim \mathbf{P}_{(j,a_j)}\) and send it;\;
                }   \ElseIf{"Productive step"} {
                    \(k = k+1\);\;
                    \(c^{k}_{(j, a_j)} = c^{k-1}_{(j,a_j)} + \eta\);\;
                } \ElseIf{"Non-productive step"} {
                    Receive \(s\) and \(i\);\;
                    \(k = k+1\);\;
                    \(c^k_{(j,a_j)} = c^{k-1}_{(j,a_j)} - \eta(1 + \widetilde{\bP}_{(j,a_j),s } - \widetilde{\bP}_{(j,a_j), i} + \mathbb{I}\{j = i\} - \mathbb{I}\{j = s\})\);
                }
            }
      }
\end{algorithm}

\textit{Remark 1.} Notice that complexity of preprocessing and Algorithm \ref{alg:amdp_solving} matches up to logarithmic factors.

\textit{Remark 2.} Messages "Productive step" and "Non-productive step" ensure that nodes update their values to actual ones. 

\textit{Remark 3.} The communication costs on each round of communication are low: each text message could be send using \(O(1)\) bits, and each message with a sample could be sent using only \(O(\log \vert \St \vert)\) bits.

\textit{Remark 4.} Additional advantage of the algorithm is a sparsity of updates: there are at most 2 values in the vector \(h\) updated each iteration. From the point of view of external memory algorithms, it gives us a possibility to make only 2 requests to the memory if the state-space is too large to store a vector \(h\) in RAM.

\textit{Remark 5.} We use a very coarse bound on \(\overline{\Theta}^2 \sim \vert \St \vert\). In practice we might expect that \(\overline{\Theta}^2 \sim \poly\log \vert \St \vert\) and, thus, reduce our dependence on \(\vert \St \vert\) to logarithmic.

\textit{Proof.}
    % First of all, we ensure that each \((j,a_j)\) worker nodes computes \(c^k = g^{(j,a_j)}_\delta(\bar v^k, h^k)\) correctly. Induction on \(k\), the basis \(c^0 = g^{(j, a_j)}_\delta(0,0) = r_{(j,a_j)}\) is correct. Now we are going to provide the induction step. In the case of productive steps
    % \begin{align*}
    %     c^k &= c^{k-1} + \eps/64 \\
    %     &= r_{(j,a_j)} - \bar v^{k-1} + (\langle \widetilde \bP_{(j,a_j)}, h^{k-1} \rangle - h^{k-1}_j) + \eps/64.
    % \end{align*}
    % It is correct since \(v^k = v^{k-1} - \eps/64\) and \(h^{k} = h^{k-1}\). In the case of non-productive steps
    % \begin{align*}
    %     c^k &= c^{k-1} - \eps/64(1 + (\widetilde{\bP}_{(j,a_j),s } -\widetilde{\bP}_{(j,a_j), i}) + 1 ) \\
    %     &= r_{(j,a_j)} - (\bar v^{k-1} + \eps/64) \\
    %     &\ \ \ + \langle \widetilde{\bP}_{(j,a_j)}, h^{k-1} - \eps/64(\eps_s - \eps_i)\rangle \\
    %     &\ \ \ - (h^{k-1}_j + \eps/64(\mathbb{I}\{i = j\} - \mathbb{I}\{s = j\})).
    % \end{align*}
    % Since \(v^k = v^{k-1} + \eps/64\) and \(h^{k} = h^{k-1} - \eps/64(e_s - e_i)\), the computation of constraints is correct.
    By Proposition \ref{lm:rounding}, we can produce \(\eps\)-optimal policy with probability at least \(1-\sigma/2\) by running Algorithm \ref{alg:md} on problem \eqref{eq:amdp_reform} with an accuracy \(\eps' = \eps/8\) because \(4(\eps_\Delta + \eps_g) = 4(\eps' + 2\delta) = 4\eps' + 8 \delta = 8 \eps'\). Performing union bound for probability of failure for a preprocessing and an algorithm, we have required probability of success of whole scheme.
    
    Since from the point of view of the head node Algorithm \ref{alg:amdp_solving} is essentially Algorithm \ref{alg:md}, we have needed guarantees on number of Mirror Descent iterations by computed constants \(M = O(1), \overline{\Theta}_0^2 = O(\tmix^2 \vert \St \vert), \kappa(E^*) = 1\) and Theorem \ref{th:primal-dual_convergence}. Notice that the parallel depth of the algorithm is equal to \(1\).
    
    Total running time complexity forms by additional \(\Atot\) computations on constraints on each iterations and each computation spends \(O(1)\) time by using previous values of computed constraints. The last observation connected to the fact that update of \(\bar v^k\) and \(\bar h^k\) also spends \(O(1)\) time by sparsity.
\hfill$ \square$

\section{\uppercase{Numerical Experiments}}\label{sec:exp}

In this section we perform numerical comparison of Algorithm \ref{alg:amdp_solving} in \textit{sequential setting} with algorithm described in paper \citep{jin2020efficiently}.

At first, we highlight technical features. In both mentioned algorithms, the value of mixing time \(\tmix\) is needed. However, computation of maximum over policies of discrete-valued function seems to be a very hard problem to be computed precisely. We replace the value of \(\tmix\) by its estimate \(\widehat \tmix\) computed over \(1000\) randomly generated policies. Additionally, to make precise comparison we compute optimal average-reward value \(v^*\) using LP-representation of the problem \eqref{eq:amdp_lp_vanila} for known model, and the average-reward value \(v^\pi\)  of the given policy \(\pi\) by computing a stationary distribution using known model.

For our comparison, we apply algorithms to two different environments: RiverSwim \citep{strehl2008analysis}, Access-Control Queuing task \citep{sutton1998rl}. We choose these environments because we have an exact model for them, and they require non-trivial exploration.

\paragraph{RiverSwim.}

We start with the environment description. RiverSwim is an environment with 6 states in a row with 2 actions for each state: swim to the left or swim to the right. Swimming to the left is always successful. For the first state, this action returns to itself and gives \(0.005\) reward. Swimming to the right move the agent to the right state with probability \(0.35\), make the agent stay in the current state with probability \(0.6\), and move to the left with probability \(0.05\). For the last state moving to the right state returns the agent to this last state and gives \(1\) reward. Thus, the optimal average-reward policy is always swimming to the right. 

\begin{figure}
\centering
\includegraphics[width=0.44\textwidth]{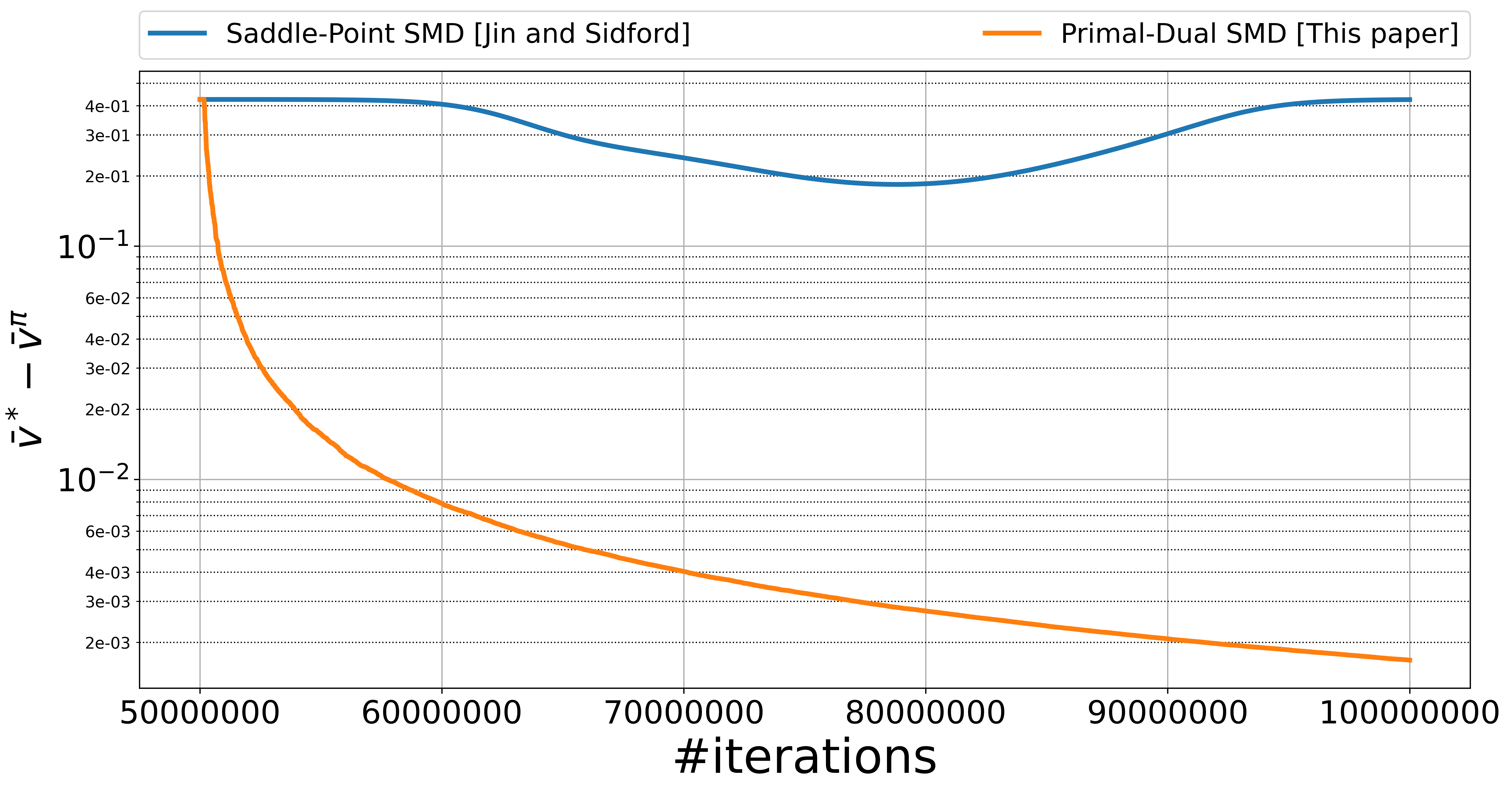}
\caption{Comparison between Stochastic Mirror Descent \citep{jin2020efficiently} and Algorithm \ref{alg:amdp_solving}  on RiverSwim environment with 6 states and 2 actions.}
\label{fig:riverswim}
 \end{figure}

Estimate value of \(\widehat{\tmix}\) is equal to \(155\) and it makes this environment hard for our algorithm \ref{alg:amdp_solving} and Stochastic Mirror Descent (SMD) \citep{jin2020efficiently}. We compare these algorithm with \(\eps = 10^{-2}\) and use \(1000\) samples for preprocessing, the result is presented on Figure \ref{fig:riverswim}. \new{Significant superior of our approach could be explained by independence of stepsize of Algorithm ~\ref{alg:amdp_solving} to the value of $t_{\text{mix}}$, whereas stepsize in algorithm of \citep{jin2020efficiently} depends on this quanitity as $1/t_{\text{mix}}^2$. In the case of this MDP, stepsize of SMD is smaller in $2.5 \cdot 10^4$ times. Additionally, notice that we generated a very little number of samples during preprocessing, and therefore these computation does not affect the total running time. Additionally, high sparsity of rewards explains the high value of the objective function during first iterates.}

%However, it requires many iterations to make the value of \(v^* - v^{\pi}\) less than \(10^{-1}\). It connected to high sparsity of rewards of this MDP. 
%Talking about QVI algorithm, for the model computed on \(1000\) samples per each state-action pair, it finds an optimal AMDP-value in a few iterations.

\paragraph{Access-Control Queuing task.} 

This environment is modelling very practical situation for using average-reward MDP model. For the description we refer to \citep{sutton1998rl}. The only difference is normalization of rewards to make them in the interval \([0,1]\).

\begin{figure}
\centering
\includegraphics[width=0.44\textwidth]{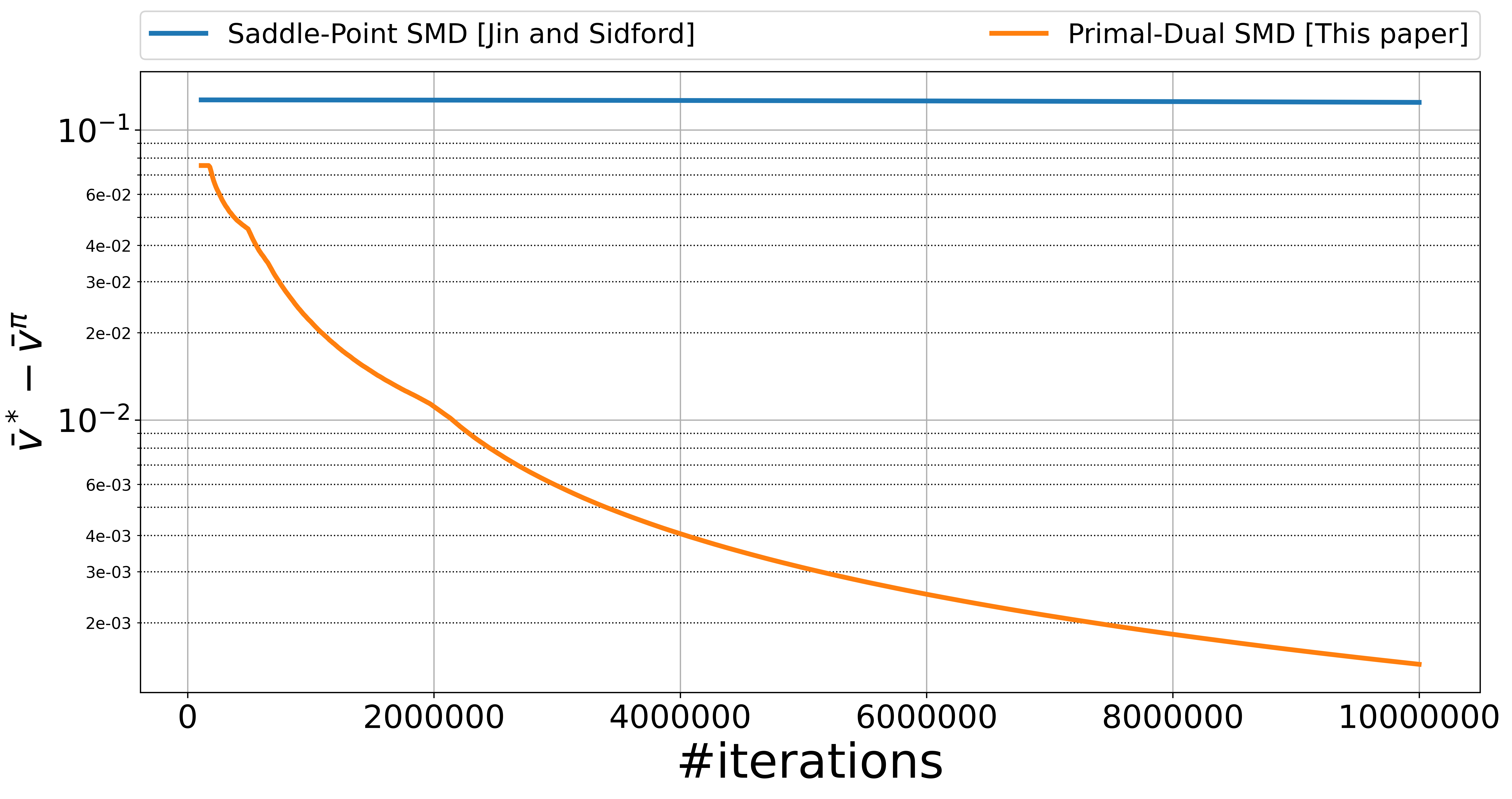}
\caption{Comparison between Stochastic Mirror Descent \citep{jin2020efficiently} and Algorithm \ref{alg:amdp_solving} on Access-Control Queuing task with 10 servers.}
\label{fig:access_control}
 \end{figure}

Estimated value of mixing time \(\widehat{\tmix} = 44\) and has dense rewards. On Figure \ref{fig:access_control} we present comparison between Algorithm \ref{alg:amdp_solving} and SMD with \(\eps=10^{-2}\) and \(500\) samples for preprocessing. \new{Again, since the value of $\tmix$ is relatively large for this MDP, worse convergence of SMD is explained by smaller stepsize that is required by theoretical analysis. Our algorithm is almost agnostic to the value of $\tmix$ and it makes it much more effective for large values of $\tmix$.}
%For this environment QVI algorithm has a very fast convergence, as in the previous case: it learned an optimal policy using \(500\) samples in a few iterations.

% \paragraph{Worst-case MDP.} This environment is the worst-case environment for the average-reward problem and it is described in \citep{jin2021tight}. This environment is natural to compare Algorithm \ref{alg:amdp_solving} and QVI because of its hardness of the second algorithm: from the theoretical point of view it has dependence \( \eps^{-3}\) on the suboptimality gap \(v^* - v^\pi\). 

% For this experiment, we used a theoretical value of mixing time \(\tmix = \lceil \log(4)/(1-\gamma) \rceil\) and environment parameters \(N=20\), \(K=5\), \(\gamma = 0.7\), and \(\eps = 10^{-3}\). Additionally, we use an additional heuristic for Algorithm \ref{alg:amdp_solving}: we memorize all samples during the running time to make computation of constraints more exact. The result is presented in Figure \ref{fig:worst_case}.

% \begin{figure}[ht!]
% \centering
% \includegraphics[width=0.44\textwidth]{images/worst_case.png}
% \caption{Comparison between Algorithm \ref{alg:amdp_solving} and Stochastic Mirror Descent \citep{jin2020efficiently} on worst-case MDP.}
% \label{fig:worst_case}
%  \end{figure} 

\section{\uppercase{Conclusion}}\label{sec:conclusion}

In this work, we proposed a parallel algorithm for solving an average-reward MDP. As far as we know, it is the first parallel algorithm in the generative model setting without reduction to the discounted problem. The interesting properties of the provided method are very low communication costs between the head node and all other nodes and the sparsity of updates that offer a possibility to work with a very large state space. 

Another contribution is the development of Mirror Descent with constraints algorithms. We provide an algorithm that works with inexact computation of constraints and prove its primal-dual properties. The setting of inexact computation of constraints was developed in \citep{lan2020algorithms} but results on primal-dual convergence of such algorithms appeared in known literature only in a deterministic exact case \citep{bayandina2018mirror}.

Turning to possible extensions, there arise natural questions.

Could a preprocessing step be avoided and make the algorithm model-free? In the current version, we needed to do a required number of preprocessing iterations to guarantee the condition on constraints. It seems possible to use an unbiased stochastic oracle for constraint evaluation. 

Another question is connected to the total work complexity. The cost of a high level of parallelism is worse total running time and sample complexity in comparison to algorithms based on saddle-point formulations. \new{It is connected to two theoretical issues: 1) $\ell_1$-approximation of the model and 2) performing Mirror Descent with box constraints. The first issue could be resolved using another approximation metric that respects MDP structure as it was done in \citep{agarwal2020model}. The second issue is fundamental for non-linear optimization due to existing lower bounds \citep{guzman2015onlc}. The only way to avoid it is to use the linear structure of mixing AMDP problem.} Could the total work time and sample complexity be reduced without losing the possibility to run in parallel?

\new{
The last question is about further generalizations of our approach. One of the interesting directions is to use our algorithm for solving constrained MDP as it is done in the paper \citep{jin2020efficiently}. However, we note that it is possible to use non-linear constraints in our case.
}

%The last question is about the choice of step-sizes. It seems possible to use an adaptive step-sizes scheme as in \citep{bayandina2018mirror,stonyakin2019adaptive}, and it might increase the practical value of the algorithm.

\subsubsection*{Acknowledgements}
The work of D.~Tiapkin was supported by the grant for research centers in the field of AI provided by the Analytical Center for the Government of the Russian Federation (ACRF) in accordance with the agreement on the provision of subsidies (identifier of the agreement 000000D730321P5Q0002) and the agreement with HSE University  No. 70-2021-00139.

\bibliography{ref}

\input{supplement}

\end{document}

%% file: supplement.tex
%%%%%%%%%%%%%%%%%%%%%%%%%%%%%%%%%%%
%%%%%% SUPPLEMENT (OPTIONAL) %%%%%%
%%%%%%%%%%%%%%%%%%%%%%%%%%%%%%%%%%%

\clearpage
\appendix

\thispagestyle{empty}

% For one-column format, uncomment the following:
\onecolumn \makesupplementtitle

\section{MISSING PROOFS}

% The supplementary materials contain detailed proofs of the results that are missing in the main paper.

\subsection{Proof of Theorem 1}
    To simplify the notation, we introduce the stochastic gradient noise function
    \begin{align}\label{eq:stoch_grad_noise_func}
        \gamma_k(y) = \begin{cases}
            \eta \langle  \hat \nabla_k f - \nabla_k f, y-x^k\rangle, & k \in I; \\
            \eta \langle \hat \nabla_k g^{(l(k))} - \nabla_k g^{(l(k))}, y-x^k\rangle, & k \in J.
        \end{cases}
    \end{align}
    The main property of this quantity is that it forms a martingale-difference sequence.
    
    Now we are going to provide some useful properties of Algorithm 1 in terms of bounds of error in objective function \(f\) and constraint satisfiability. Afterwards, we combine all properties to a convergence analysis of the algorithm. But at first, we state one useful technical result.
    We will actively use the following lemma:
    \begin{lemma}[\cite{bayandina2018mirror}]\label{lm:md_main}
        Let \(h\) be some convex function over a set \(Q\), \(\eta > 0\) is a stepsize, \(x \in Q\). Let the point \(x^+ = \Mirr(x, \eta(\nabla h(x) + \Delta))\), where \(\Delta\) is some vector from the dual space. Then, for any \(y \in Q\)
        \begin{align*}
            \eta(h(x) - h(y) + \langle \Delta, x-y \rangle) \leq \eta \langle \nabla h(x) + \Delta, x - y \rangle\\
            \leq \frac{\eta^2}{2} \Vert \nabla h(x)+ \Delta \Vert^2_* + V(x,y) - V(x^+, y).
        \end{align*}
    \end{lemma}
\begin{lemma}\label{lm:f_basic_bound}
        For a point \(\hat x\) produced by Algorithm 1 and any \(y \in Q\) the following holds
        \begin{align}\label{eq:f_basic_bound}
        \begin{split}
                \eta \vert I \vert \cdot  (f(\hat x) - f(y)) &\leq \frac{\eta^2 M^2}{2} N + [d(y) - d(x^0)] - \vert J \vert \eta \cdot \eps \\
                &+ \sum_{k=0}^{N-1} \gamma_k(y) + \eta \sum_{k \in J} g^{(l(k))}(y).\\
        \end{split}
        \end{align}
    \end{lemma}
    
    \textit{Proof.}
            By the construction of "productive" and "non-productive" steps in Algorithm 1 and Lemma \ref{lm:md_main}, we have for all \(y \in Q\)
    \begin{align*}
        \begin{split}
            \eta(f(x^k) - f(y)) &\leq \frac{\eta^2 M^2 }{2}+ V(x^k, y) - V(x^{k+1}, y) \\
            &+  \eta \langle \hat \nabla_k f  - \nabla_k f, y-x^k\rangle, 
        \end{split} & k \in I; \\
        \begin{split}
            \eta(g^{(l(k))}(x^k) - g^{(l(k))}(y)) &\leq \frac{\eta^2 M^2}{2} + V(x^k, y) - V(x^{k+1}, y) \\
            &+  \eta \langle \hat \nabla_k g^{(l(k))} - \nabla_k g^{(l(k))}, y-x^k\rangle
        \end{split} & k \in J.
    \end{align*}

        By definition of \(\delta\)-approximation, we have the following inequalities for "productive" and "non-productive" steps respectively
    \begin{align*}
         \eta(f(x^k) - f(y)) &\leq \frac{\eta^2 M^2}{2}  + V(x^k, y) - V(x^{k+1}, y) + \gamma_k(y),\\
         \eta( g_\delta^{(l(k))}(x^k) - g^{(l(k))}(y)) &\leq \frac{\eta^2 M^2}{2} + V(x^k, y) - V(x^{k+1}, y) + \gamma_k(y)
         + \eta \delta .
    \end{align*}
    
    Sum all these inequalities over all \(k \in I\) and \(k \in J\) and use the fact that \(I \cup J = \{0,\ldots,N-1\}\)
    \begin{align}\label{eq:f_prelim_bound}
        \begin{split}
            &\sum_{k \in I} \eta(f(x^k) - f(y)) + \sum_{k \in J} \eta(g_\delta^{(l(k))}(x^k) - g^{(l(k))}(y)) \\
        &\leq \frac{\eta^2 M^2}{2} \vert I \vert  +  \frac{\eta^2 M^2}{2} \vert J \vert + \sum_{k=0}^{N-1} [V(x^k,y) - V(x^{k+1},y)] + \sum_{k=0}^{N-1} \gamma_k(y) + \vert J \vert \eta  \delta.   
        \end{split}
    \end{align}
    By the choice of \(x^0 = \argmin_{x \in Q} d(x)\), we have \[
        \sum_{k=0}^{N-1} [V(x^k,y) - V(x^{k+1},y)] \leq  V(x^0,y) = d(y) - d(x^0) - \langle \nabla d(x^0), y - x^0\rangle = d(y) - d(x^0).
    \]
    Using definition of "non-productive" steps \(g_\delta^{(l(k))}(x^k) > \eps + \delta\) and convexity of \(f\)
    \begin{align*}
        \sum_{k \in I} \eta(f(x^k) - f(y)) + \sum_{k \in J} \eta( g_\delta^{(l(k))}(x^k) - g^{(l(k))}(y)) \\
        \geq \eta \vert I \vert (f(\hat x) - f(y)) + \eta \vert J \vert (\eps + \delta) - \eta \sum_{k \in J} g^{(l(k))}(y).
    \end{align*}
    By application of inequality \eqref{eq:f_prelim_bound} and simple regrouping of terms, we finish the proof.
\hfill$ \square$

During this section we set \(y = x^*\). In this case \(d(x^*) - d(x^0) \leq \Theta_0^2\). Further we are going to derive bound on \(g^{(l)}\) at \(\hat x\) for each separate function \(g^{(l)}\):
\begin{lemma}\label{lm:g_basic_bound}
    For \(\hat x\) produced by Algorithm 1 the following holds
    \begin{equation}\label{eq:g_basic_bound}
         g^{(l)}(\hat x) \leq \eps + 2\delta.
    \end{equation}
\end{lemma}
\textit{Proof.}
    By convexity and definition of a \(\delta\)-approximation
    \[
    g^{(l)}(\hat x)  \leq \frac{1}{\vert I\vert} \sum_{k \in I} g^{(l)}(x^k) \leq \frac{1}{\vert I\vert} \sum_{k \in I} (g^{(l)}_\delta(x^k) + \delta).
    \]
    We finish the proof by definition of "productive" steps and a set \(I\).
\hfill$ \square$

% Here we see the main difficulty: we have to control the mean of \(\hat \gamma_k^{(l)}\) over a random set of indices to control \(g^{(l)}(\hat x)\). This sum is expected to be small if and only of \(\vert I\vert\) is big enough: if \(\vert I \vert = 1\), the only bound that we can produce is \(\eps + 2 \delta + B\) that can be very large. Thus, we have to guarantee that \(\vert I\vert\) will be at least linear on \(N\) with high probability. It will be the moment where the technical assumption \textbf{(A4')} appears.

Before obtaining final bounds on \(f\) we prove some technical fact:
\begin{lemma}\label{lm:norm_bound_theta}
    For any \(y \in Q\): \(\Vert y - x^0 \Vert^2 \leq 2(d(y) - d(x^0))\).
\end{lemma}
\textit{Proof.}
    Follows from strongly convexity of  \(d\) with respect to norm \(\Vert \cdot \Vert\) and the fact that \(x^0\) is a minimum of \(d\).
\hfill$ \square$

Now we derive bounds on sums of our stochastic noise function with high probability. 
\begin{lemma}\label{lm:azuma_primal}
    Define event \(\mathcal{E}\) such that the following inequalities holds
    \begin{align*}
        \sum_{k=0}^{N-1} &\gamma_k(x^*) < \frac{2\Theta_0}{M} \cdot \sqrt{2N \eps^2 \log(1/\sigma)}.
    \end{align*}
    Then \(\Prob[\mathcal{E}] \geq 1 - \sigma\).
\end{lemma}
\textit{Proof.}
    Define filtration of \(\sigma\)-algebras \(\mathcal{F}_k = \sigma\left(\{ \xi^i, \xi^i_{(l)}\}_{i \leq k}\right)\). Notice that \(\gamma_k(x^*)\) is a martingale-difference sequence adapted to \(\mathcal{F}_k\).
    
    Let us derive bound on the sum. Notice that by Holder inequality and Lemma \ref{lm:norm_bound_theta}
    \begin{align*}
        \vert \gamma_k(x^*) \vert \leq  2\eta M \Vert x^* - x^k \Vert \leq \frac{4\eps \Theta_0}{M},
    \end{align*}
    Then by Azuma-Hoeffding inequality
    \begin{align*}
        \Prob\left[\sum_{k=0}^{N-1} \gamma_k(x^*) \geq t_1 \right] &\leq \exp\left( \frac{-t_1^2}{2N \cdot (16 \eps^2 \Theta_0^2)/(M^2)} \right);
    \end{align*}
    By setting \(t_1 = 4\Theta_0 M^{-1} \cdot \sqrt{2N \eps^2 \log(1/\sigma)} \) we finish the proof.
\hfill$ \square$

One could apply Lemma \ref{lm:azuma_primal} to inequalities in Lemmas \ref{lm:f_basic_bound} and \ref{lm:g_basic_bound}:
\begin{corollary}\label{cor:final_bounds}
    Under event \(\mathcal{E}\) defined in Lemma \ref{lm:azuma_primal}, the following inequalities holds for all \(l \in [m]\)
    \begin{align*}
        \eta \vert I \vert (f(\hat x) - f(x^*)) &< \eta \vert I\vert \eps  - \frac{\eps^2 N}{2 M^2}  + \Theta_0^2  +  \frac{4\Theta_0  \sqrt{2N\eps^2 \log(1/\sigma)}}{M};\\
        g^{(l)}(\hat x) &< \eps + 2\delta.
    \end{align*}
\end{corollary}
\textit{Proof.}
    The second inequality follows directly from Lemma \ref{lm:azuma_primal}. The first one uses definition \(\eta = \eps/M^2\), inequalities \(g^{(l))}(x^*) \leq 0\) and, finally, Lemma \ref{lm:azuma_primal}.
\hfill$ \square$

Now we are going to prove last technical result before the statement of the main theorem of this subsection:
\begin{lemma}\label{lm:quadratic_solve_primal}
    If \(\sigma < 1\) and
    \[
        N \geq N_0 = \frac{280 \cdot \Theta_0^2 M^2 \log(1/\sigma)}{\eps^2} 
    \]
    then
    \[
        \frac{\eps^2 N}{2 M^2} - \Theta_0^2 - \sqrt{2N\eps^2}\left( \frac{4\Theta_0 \sqrt{\log(1/\sigma)}}{M} \right) \geq 0.
    \]
\end{lemma}
\textit{Proof.}
    Fix variable \(t = 2 N \eps^2\). Then we have quadratic inequality of type \(at^2 - (b_1 + b_2)t - c \geq 0\). By exact formula for quadrativ equation and Taylor approximation we know that inequality holds if \(t > (b_1 + b_2)/a + c/(b_1 + b_2)\). By numeric inequality \(1/(b_1+b_2) \leq 1/b_1\) for nonnegative \(b_1,b_2\), we have that enough to choose \(t \geq b_1/a + b_2/a + c/b_1\). Using our choice of \(N\) and numeric inequality \(2a^2 + 2b^2 \geq (a+b)^2\)
    \begin{align*}
        \sqrt{2N \eps^2} &\geq \sqrt{280} \cdot \Theta_0 M \sqrt{\log(1/\sigma)} \geq 16 \Theta_0 \sqrt{\log(3/\sigma)} \cdot M + \frac{\Theta_0^2 M}{2 \Theta_0 \sqrt{\log(1/\sigma)}}.
    \end{align*}
    Here we use that \(\log(1/\sigma) \geq 1\) and \(16 + 1/2 \leq \sqrt{280}\).
\hfill$ \square$

Finally, we are ready to state theorem that describe convergence properties of Algorithm 1.

\textit{Proof of Theorem 1.}
    To guarantee that \(f(\hat x) - f(x^*) \leq \eps\) with probability at least \(1-\sigma\), we use the first inequality in Corollary \ref{cor:final_bounds} and Lemma \ref{lm:quadratic_solve_primal}. To guarantee satisfaction of constraints, we simply use Corollary \ref{cor:final_bounds}.
\hfill$ \square$

\subsection{Proof of Theorem 2}

Recall our estimate of a dual variables
\begin{equation*}%\label{eq:dual_variables}
    \hat \lambda_l = \frac{1}{\vert I \vert }\sum_{k \in J}  I\{l = l(k)\}.
\end{equation*} 

\begin{lemma}\label{lm:gap_basic_bound}
    Suppose \(\hat x\) is an output of Algorithm 1 and \(\hat \lambda\) is defined as in \eqref{eq:dual_variables} and \(\overline{\Theta}_0^2 \geq \sup_{y \in Q} (d(y) - d(x^0))\).
    
    Then the following holds
    \begin{align}\label{eq:gap_basic_bound}
        \begin{split}
                \eta \vert I \vert \Delta(\hat x, \hat \lambda)  &\leq \frac{\eta^2 M^2}{2} N  + \overline{\Theta}_0^2 - \vert J \vert \eta \eps \\
                &+ \sum_{k=0}^{N-1} \gamma_k(x^0) + \overline{\Theta}_0 \sqrt{2} \cdot \left\Vert \sum_{k=0}^{N-1} \Delta_k \right\Vert_*,
        \end{split}
    \end{align}
    where \(\Delta_k\) is defined as follows
    \[
      \Delta_k = \begin{cases}
            \eta \left(  \hat \nabla_k f - \nabla_k f \right), & k \in I \\
            \eta \left(  \hat \nabla_k g^{(l(k))} - \nabla_k g^{(l(k))} \right), & k \in J.
        \end{cases}
    \]
\end{lemma}
\textit{Proof.}
    We start from Lemma \ref{lm:f_basic_bound}. Here we move all terms consist of \(y\) to the right-hand side and minimize over \(y\)
    \begin{align*}
        \eta \vert I \vert f(\hat x) &\leq \frac{\eta^2 M^2}{2} N - \vert J \vert \eta \eps \\
        &+ \min_{y \in Q} \left\{ d(y) - d(x^0) + \sum_{k=0}^{N-1} \gamma_k(y) + \eta \vert I \vert f(y) + \eta \sum_{k \in J} g^{(l(k))}(y) \right\}.
    \end{align*}
    
    Notice that by definition of \(\hat \lambda\) we have
    \[
        \eta \vert I \vert \phi(\hat \lambda) = \min_{y \in Q} \left\{  \eta \vert I \vert f(y) + \eta \sum_{k \in J} g^{(l(k))} \right\}.
    \]
    Thus, we have to upper bound \(d(y) - d(x^0)\) and \(\sum_{k=0}^{N-1} \gamma_k(y)\) without dependence on \(y\) to obtain required result. The first upper bound is trivial: \(d(y) - d(x^0) \leq \overline{\Theta_0}^2\) by definition of \(\overline{\Theta_0}^2\).
    
    To analyse the second term, we use the definition of \(\gamma_k\) in terms of \(\Delta_k\) and Holder inequality
    \[
        \sum_{k=0}^{N-1} \gamma_k(y) = \sum_{k=0}^{N-1} \langle \Delta_k, y - x^0 + x^0 -  x^k \rangle \leq \Vert  y - x^0 \Vert \left\Vert \sum_{k=0}^{N-1} \Delta_k \right\Vert_* + \sum_{k=0}^{N-1} \gamma_k(x^0).
    \]
    The last step is to apply Lemma \ref{lm:norm_bound_theta} and definition of \(\overline{\Theta}_0^2\) to obtain uniform bound on \(\Vert y - x^0 \Vert\).
\hfill$ \square$

Our next goal is to derive bound on the right-hand side of \eqref{eq:gap_basic_bound}. It is possible using concentration of measure techniques as in Lemma \ref{lm:azuma_primal}.

\begin{lemma}\label{lm:azuma_dual}
    Define event \(\mathcal{E}'\) such that the following inequalities holds
    \begin{align*}
        \sum_{k=0}^{N-1} \gamma_k(x^0) &< \frac{2\overline{\Theta}_0}{M} \cdot \sqrt{2N \eps^2 \log(2/\sigma)}\\
        \left\Vert \sum_{k=0}^{N-1} \Delta_k\right\Vert_* &< \frac{\sqrt{2 \kappa(E^*) } + \sqrt{4\log{(2/\sigma)}}}{M}\sqrt{2N\eps^2} ,
    \end{align*}
    where \(\Delta_k\) is defined in \ref{lm:gap_basic_bound} and \(\kappa(E^*)\) is a constant of Nemirovski's inequality \citep{boucheron2013concentration} for the dual space. Then \(\Prob[\mathcal{E}'] \geq 1 - \sigma\).
\end{lemma}
\textit{Remark.} If \(E\) has a finite dimension \(d\), then we always have \(\kappa(E^*) \leq d\). Additionally, if \(E\) is endowed with \(\ell_p\) norm, then \(E^*\) is endowed with \(\ell_q\) norm, where \(1/p + 1/q = 1\), and there is a more precise bound, according to \citep{lutz2010nemirovski}
\[
    \kappa(E^*) \leq K\left(\frac{p}{p-1}, d\right) = \begin{cases}
        d^{\frac{2}{p}-1}, & p \in [1,2]\\
        d^{1-\frac{2}{p}}, & p \in (2,+\infty] \\
    \end{cases}
\]
In particular, if \(E\) has \(\ell_2\) norm, \(\kappa(E^*) = 1\). For \(p\in[2,+\infty]\) this bound is tight, however, in the case \(p\in[1,2]\) and \(d\geq 3\) it could be improved to, for instance, a logarithmic bound \(\kappa(E^*) \leq 2e \log(d)- e\), that could be useful in the case of \(\ell_1\)-norm and an entropy prox-function.
\textit{Proof.}
    The case of first inequality is identical to Lemma \ref{lm:azuma_primal} by rescaling \(\sigma\) to \(\sigma/2\).
    
    To ensure the last inequality, we apply bounded difference inequality \citep{boucheron2013concentration}. This follows by observing that \(Z = \Vert \sum_{k=0}^{N-1} \Delta_k \Vert_*\) satisfies bounded difference condition with a constant \(4 \cdot \eps \underline{M}^{-1}\) that is greater than \( 2 \cdot \Vert \Delta_k \Vert_*\) a.s.. Thus
    \[
        \Prob[Z - \E Z > t_2]  \leq \exp\left(\frac{-t_2^2}{2\cdot N \eps^2 \cdot 4 M^{-2}} \right).
    \]
    Take \(t_2 = 2M^{-1} \sqrt{2 N \eps^2} \cdot \sqrt{\log(4/\sigma)}\) and the last we have to do is to bound expectation of \(Z\). Here we apply Nemirovski's inequality
    \begin{align*}
        \left( \E \left\Vert \sum_{k=0}^{N-1} \Delta_k \right\Vert_* \right)^2 \leq \E\left[\left\Vert \sum_{k=0}^{N-1} \Delta_k \right\Vert_*^2\right] \leq \kappa(E^*) \sum_{k=0}^{N-1} \E\left[ \left\Vert  \Delta_k \right\Vert_*^2\right] \leq  \kappa(E^*) \frac{4 N \eps^2}{M^2}.
\end{align*}
By application of the union bound we finish the proof.
\hfill$ \square$

\begin{corollary}\label{cor:gap_final_bound}
    Under event \(\mathcal{E}'\) defined in Lemma \ref{lm:azuma_dual}, the following inequalities holds for all \(l \in [m]\)
    \begin{align*}
        \eta \vert I \vert \Delta(\hat x, \hat \lambda) &< \eta \vert I\vert \eps  - \frac{\eps^2 N}{2 M^2} + \overline{\Theta}_0^2 + \sqrt{2N\eps^2}\left( \frac{2\overline{\Theta}_0( 4\sqrt{\log(2/\sigma)} + \sqrt{2\kappa(E^*)} )}{M}\right);\\
        g^{(l)}(\hat x) &< \eps + 2\delta.
    \end{align*}
\end{corollary}
\textit{Proof.}
    Inequalities on \(g^{(l)}(\hat x)\) is equivalent to the same in \ref{cor:final_bounds}. Inequality on \(\Delta(\hat x, \hat \lambda)\) follows from a combination of Lemma \ref{lm:gap_basic_bound} and Lemma \ref{lm:azuma_dual}.
\hfill$ \square$

Now we are ready to state a technical lemma that is similar to Lemma \ref{lm:quadratic_solve_primal}.

\begin{lemma}\label{lm:quadratic_solve_dual}
    If \(\sigma \leq 1/2\) and 
    \[
        N \geq N_0' = \frac{128 \overline{\Theta}_0^2 M^2 (17 \log(2/\sigma) + 2 \kappa(E^*))}{\eps^2}
    \]
    then
    \[
        \frac{\eps^2 N}{2 M^2} - \overline{\Theta}_0^2 - \sqrt{2N\eps^2}\left( \frac{2\overline{\Theta}_0( 4\sqrt{\log(2/\sigma)} + \sqrt{2\kappa(E^*)} )}{M} \right) \geq 0.
    \]
\end{lemma}
\textit{Proof.}
    Using the same reasoning about solving quadratic inequality in terms of \(t = \sqrt{2 N\eps^2}\) as in Lemma \ref{lm:quadratic_solve_primal} it is sufficient to show that
    \begin{align*}
        \sqrt{2N \eps^2} &\geq 4 \left(2 \overline{\Theta}_0 M (4 \sqrt{\log(2/\sigma)} + \sqrt{2\kappa(E^*)})\right) \\
        &+ \frac{\overline{\Theta}_0^2 M}{2 \overline{\Theta}_0(4 \sqrt{\log(2/\sigma)} + \sqrt{2\kappa(E^*))}}.
    \end{align*}
    Since \(\sigma \leq 1/2 \Rightarrow \log(2/\sigma) \geq 1\) and \(\sqrt{2\kappa(E^*)} \geq 1\), we have that it is sufficient to show that
    \[
        \sqrt{2N \eps^2} \geq 4 \left(2 \overline{\Theta}_0 M ((4+1/20) \sqrt{\log(2/\sigma)} + \sqrt{2\kappa(E^*)} ) \right).
    \]
    By numeric inequality \(2a^2 + 2b^2 \geq (a+b)^2\) and \((4+1/20)^2 \leq 17\), it is satisfied with our choice of \(N \geq N_0'\).
\hfill$ \square$

Now we are ready to prove our main result.

\textit{Proof of Theorem 2.}
    The inequality on \(g(\hat x)\) is satisfies by Corollary \ref{cor:gap_final_bound}. We have to satisfy inequality on duality gap. The inequality on duality gap follows directly from Corollary \ref{cor:gap_final_bound} and Lemma \ref{lm:quadratic_solve_dual} since \(N \geq N_0'\). Notice that a probability \(1-\sigma\) appears from the event \(\mathcal{E}'\) defined in Lemma \ref{lm:azuma_dual}.
\hfill$ \square$

\subsection{Proof of Proposition 1}
At first, we prove one technical lemma. 
\begin{lemma}[Estimation of parameters of categorical distribution]
    Suppose that we have samples \(\{X_i\}_{i=1}^N, X_i \in \{1,\ldots,d\}\) drawn from categorical distribution with a parameter \(s \in \Delta^d\). Define an empirical estimate \(\hat s = \frac{1}{N} \sum_{i=1}^N e_{X_i} \in \Delta^d \), where \(e_j\) is an \(j\)-th standard basis vector. 
    
    Then, for any \(\delta',\sigma' > 0\), the inequality \(\Vert s - \hat s \Vert_1 \leq \delta'\) holds with probability at least \(1-\sigma'\) if 
    \[
        N \geq  \frac{8d + 4 \log(1/\sigma')}{\delta'^2}.
    \]
\end{lemma}
\textit{Proof.}
    First of all, notice that \(\E e_{X_i} = s\). Denote by \(\Delta\) a centred random variable \(\hat s - s = \frac{1}{N} \sum_{i=1}^n(e_{X_i} - s)\). Then, by Nemirovski's inequality and bound on the \(\ell_1\) norm of elements of a simplex, we might estimate the mean of the square of  \(\ell_1\) norm of this random variable: \( \E \Vert \Delta \Vert^2_1 \leq 4d/N \).
    
    Let us define function \(f(X_1,\ldots,X_N) = \Vert \Delta \Vert_1\). It might be checked that this function satisfies conditions of the bounded difference inequality \citep{boucheron2013concentration} with constant \(2/N\). Thus, for all \(t > 2 \sqrt{d/N} \geq \E \Vert \Delta \Vert_1\)
    \[
        \Prob[\Vert \Delta \Vert_1 > t] = \Prob[\Vert \Delta \Vert_1 - \E \Vert \Delta \Vert_1 > t  - \E \Vert \Delta \Vert_1] \leq \exp\left(\frac{-(t- 2\sqrt{d/N})^{2}}{2 N^{-1}} \right).
    \]
    Taking \(N \geq (8d + 4 \log(1/\sigma'))\cdot (\delta')^{-2}\) and \(t =\delta'\), we finish the proof.
\hfill$ \square$

Using this simple lemma, we can easily obtain required sample complexity for the preprocessing even in parallel sampling setting. Remember that we are going to use Algorithm 1, hence, we want to approximate each constraint function.

\textit{Proof of Proposition 1.}
    Notice that each \(\bP_{(a)} \in \Delta^{\vert \St \vert}\) is a parameters of a categorical distribution we sampling from, and the estimator from the previous lemma could be applied.
    
    By Holder's inequality and the definition of the search space for \(h\) we have that
    \(
        \vert \langle \bP_{(a)} - \widetilde \bP_{(a)}, h \rangle \vert \leq 2M \cdot \Vert \bP_{(a)} - \widetilde \bP_{(a)} \Vert_1.
    \)
    Hence, to make it less than \(\delta'\), we required to make \(\ell_1\) norm of difference less than \(\delta'/(2M)\). To make all conditions work simultaneously, it is enough to make the probability in the terms of the previous lemma less than \(\sigma' / \Atot\) and apply the union bound over all \(\Atot\) conditions. The last observation that finishes the proof is that sampling for an estimation of each \(\widetilde \bP_a\) could be done separately and independent.
\hfill$ \square$

\subsection{Proof of Proposition 2}
\textit{Proof.}
Conditions of the proposition give us the following guarantees in terms of the duality gap and constraint satisfaction with needed probability \(\geq 1 - \sigma\)
\begin{align}\label{eq:amdp_constraint}
    \begin{split}
    &\bar v^\eps - \min_{\bar v, h}\left( \bar v + (\mu^\eps)^\top \left( - \bar v \one - (\hat\bI- \bP)h + \br \right)\right)  \leq \eps_f, \\ 
    &\bar v^\eps \one + (\hat\bI- \bP) h^\eps - \br \geq -\eps_g \one.
    \end{split}
\end{align}

We can rewrite the first condition in more suitable terms of \(\lambda^\eps\)
\begin{align}
    \begin{split}
    \forall \bar v, \forall h \in \mathcal{X}: &\eps_f \geq (\bar v^\eps - \bar v) + \bar v \cdot  (\lambda^\eps)^\top \one + (\lambda^\eps)^\top ( [I- \bP^\pi] h - \br^\pi).
    \end{split}\label{eq:amdp_duality_gap}
\end{align}
Now we have all required instruments and we can bound the expectation of an average value of our policy
\begin{align*}
    v^\pi &= (\nu^\pi)^\top \br^\pi = (\nu^\pi - \lambda^\eps)^\top ( [\bP^\pi - I] h^\eps + \br^\pi) +  (\lambda^\eps)^\top ( [\bP^\pi - I] h^\eps + \br^\pi ).
\end{align*}
Here we used the stationary of our policy: \((\nu^\pi)^\top (\bP^\pi - I) = 0\). For simplicity, we remind that \(\mu^\eps, \lambda^\eps \geq 0\), and, hence, \(\langle \mu^\eps, \one\rangle = \langle \lambda^\eps, \one \rangle = \Vert \lambda^\eps \Vert_1\).

Firstly, bound the second term by \ref{eq:amdp_duality_gap}
\[
     (\lambda^\eps)^\top ( [\bP^\pi - I] h^\eps + \br^\pi ) \geq\bar v^\eps - \bar v(1 - \Vert \lambda^\eps\Vert_1) - \eps_f.
\]
To bound the first term we will use Lemma 7 from \citep{jin2020efficiently}:
\[
    \Vert (\bI - \bP^\pi + \one (\nu^\pi)^\top)^{-1} \Vert_{\infty } \leq M.
\]
Also we have for all \(\mu \geq 0\) by primal feasibility:
\(
     \mu^\top ([\hat \bI - \bP] h^* - \br + \bar v^* \one ) \geq 0.
\)
We can combine it with our condition \eqref{eq:amdp_duality_gap} for arbitrary \(h\) and \(\bar v = \bar v^*\), using the fact that \(\bar v^\eps \geq \bar v^* - \eps_g\)
\begin{align*}
    \eps_f + \eps_g &\geq \bar v^* - \bar v^* + (\mu^{\eps})^\top([\hat \bI - \bP] h - \br + \bar v^* \one) \\
    &\geq  (\mu^{\eps})^\top([\hat \bI - \bP] h - \br + \bar v^* \one) - (\mu^\eps)^\top ([\hat \bI - \bP] h^* - \br + \bar v^* \one )  \\
    &= (\mu^{\eps})^\top([\hat \bI - \bP] (h - h^*)) = (\lambda^{\eps})^\top [\bI - \bP^\pi] (h - h^*).
\end{align*}

Hence, we have
\begin{align*}
    2M \Vert (\lambda^{\eps})^\top [\bI - \bP^\pi] \Vert_1 &= \max_{h \in \mathbb{B}^{\vert \St \vert}_{2M}} (\lambda^{\eps})^\top [\bI - \bP^\pi] h \\
    &= \max_{h \in \mathbb{B}^{\vert \St \vert}_{2M}} (\lambda^{\eps})^\top [\bI - \bP^\pi] (h - h^*) +  (\lambda^{\eps})^\top [\bI - \bP^\pi] h^* \\
    &\leq \eps_f + \eps_g + M \Vert (\lambda^{\eps})^\top [\bI - \bP^\pi] \Vert_1.
\end{align*}
Combining with the fact that \(\nu^\pi(\bI - \bP^\pi) = 0\), we obtain\(
     \Vert (\nu^\pi - \lambda^{\eps})^\top [\bI - \bP^\pi] \Vert_1 \leq \frac{\eps_f + \eps_g}{M}. \)
By almost the same argument as in \citep{jin2020efficiently}, we also have
\(
    \vert \langle \nu^\pi - \lambda^\eps, \br^\pi \rangle \vert \leq \eps_f + \eps_g.
\)
Hence, there is a bound on the required first term:
\[
    (\nu^\pi - \lambda^\eps)^\top ( [\bP^\pi - I] h^\eps + \br^\pi) \geq -2 M \cdot \frac{(\eps_f + \eps_g)}{M} - (\eps_f + \eps_g).
\]
Overall, we obtain the required inequality by taking \(\bar v = 0\) and by feasibility of the pair \((\bar v^\eps + \eps_g, h^\eps)\)
\[
    \bar v^\pi \geq  \bar v^\eps - \bar v(1 - \Vert \lambda^\eps\Vert_1) - \eps_f - 3(\eps_f + \eps_g) \geq \bar v^* - 4(\eps_f + \eps_g).
\]
\hfill$ \square$

%\section{ADDITIONAL EXPERIMENTS}

%

%\bibliography{ref}

%\vfill

%\end{document}